\documentclass{article}
\usepackage{amsmath,amsfonts,amssymb,amscd}
\usepackage{hyperref}
\usepackage[dvips]{graphicx}
\newcommand{\PP}{{\mathbb{P}}}

\newcommand{\C}{{\mathbb{C}}}
\newcommand{\Q}{{\mathbb{Q}}}

\newcommand{\R}{{\mathbb{R}}}
\newcommand{\Z}{{\mathbb{Z}}}

\newcommand{\calO}{\mathcal{O}}

\newcommand{\calI}{{\cal I}}
\newcommand{\calJ}{{\cal J}}

\newcommand{\calL}{{\cal L}}

\newcommand{\calR}{{\cal R}}
\newcommand{\calS}{{\cal S}}

\newtheorem{theorem}{Theorem}
\newtheorem{corollary}[theorem]{Corollary}
\newtheorem{lemma}[theorem]{Lemma}
\newtheorem{proposition}[theorem]{Proposition}

\title{Adjunction conditions for 1-forms on surfaces in projective three-space}
\author{Joseph Steenbrink\\Huygens Institute for Mathematics,\\ Astrophysics
and Particle Physics\\
Radboud University Nijmegen\\
E-mail J.Steenbrink\@@science.ru.nl}
\begin{document}
\maketitle
\begin{abstract} 
We study the relation between a certain graded part of the Jacobian ring of a 
projective hypersurface and a certain graded quotient for the Hodge filtration of its 
primitive cohomology, in the case that the hypersurface has at most isolated 
singularities. We distinguish a class of singularities for which this relation 
is best possible. The main examples refer to the surface case.  
\end{abstract}
\section*{Introduction}
Let $X \subset P = \PP^{n+1}(\C)$ be a smooth hypersurface given by a homogeneous polynomial $F$ of degree $d$. 
We let 
\begin{eqnarray*}
\calS & =\C[X_0,\ldots,X_{n+1}] \\
\calS_k& = \{G\in\calS\mid G \mbox{ homogeneous of degree }k\}\\
\calJ(F) & = \mbox{ ideal in }\calS \mbox{ generated by }\partial_0F,\ldots,\partial_{n+1}F \\
\calR & = \calS/\calJ(F)
\end{eqnarray*}
Then we have isomorphisms 
$$\calR_{d(p+1)-n-2} \rightarrow H^{n-p,p}_0(X)$$
where $H_0$ denotes primitive cohomology.

We intend to investigate the relation between  $\calR_{d(p+1)-n-2}$ and the cohomology of 
$X$ in the case where $X$ has isolated singularities. More precisely we investigate for 
which singularities there is a direct relation between $\mathrm{Gr}^{n-p}_FH^n(X,\C)$ and $\calR_{d(p+1)-n-2}$.  

It appears that  $\calR_{d(p+1)-n-2}$ is closely related with the cohomology groups 
of the sheaf $\omega_X^{n-p}$ of residues of logarithmic $n-p+1$-forms on $\PP^{n+1}$, whereas 
$\mathrm{Gr}^{n-p}_FH^n(X,\C)$ is related to the $(n-p)$-th graded piece of the \emph{filtered de Rham complex} 
$\tilde{\Omega}_X^\bullet$ of $X$. 
So our question boils down to a comparison of $\omega_X^{n-p}$ (the sheaf of  \emph{Barlet forms}) and 
$\mathrm{Gr}^{n-p}_F\tilde{\Omega}_X^\bullet$.

The case $p=0$ is classical. The $n$-th graded part of the filtered de Rham complex of $X$ 
is the sheaf of meromorphic $n$-forms on $X$ which lift holomorphically to any resolution. 
For each isolated singular point $x$ of $X$ one has an ideal $\calI \subset \calO_{X,x}$ 
with the property that the residue along $X$ of a rational $n+1$-form 
$$\omega = \frac{A\Omega}{F}$$
on $\PP^{n+1}$ with a first order pole along $X$ extends holomorphically to any resolution 
if and only if $A_x \in I$ for each $x$. These conditions are called the \emph{adjunction conditions} 
and rational singularities are characterized by the fact that they do not impose adjunction conditions, 
i.e. $\calI = \calO_{X,x}$.

In this paper we study the case $p=1$. It will become clear that the main case of interest  
is the surface case. We come to a satisfactory picture in the case that $X$ is a surface with 
only a certain class of singularities, which includes the rational double points and cusps. 
We also deal with the case of surfaces in weighted projective spaces. We will give examples 
which show that the analogous class of singularities in dimensions different from two is probably empty.

Our interest in this problem was roused by a question of Remke Kloosterman. 

\section{Differentials on singular spaces} 
We recall some facts concerning sheaves of holomorphic differentials on singular spaces. 

Recall that a V-manifold is a complex analytic space which is locally isomorphic to the quotient of a 
complex ball by a finite group of biholomorphic transformations. Local models for $n$-dimensional 
V-manifolds are of the form $B^n/G$ where $B^n$ is the $n$-dimensional open unit ball and $G$ is a 
\emph{small} finite subgroup of $\mathrm{U}(n,\C)$ i.e. no element of $G$ has $1$ as an eigenvalue of multiplicity $n-1$. 

For a V-manifold $X$ one defines sheaves $\tilde{\Omega}^p_X,\ p\geq 0$ on $X$ as follows. 
Let $\Sigma$ be the singular locus of $X$. Because a V-manifold is normal, it has codimension at least two in $X$. 
Let $j:X \setminus \Sigma \to X$ denote the inclusion map, and put 
\begin{equation}
\tilde{\Omega}^p_X := j_\ast\Omega^p_{X\setminus\Sigma}.
\end{equation}
\begin{theorem}
\begin{enumerate}
\item If $X$ is a V-manifold and $\pi:\tilde{X}\to X$ a resolution of singularities, 
then $$\tilde{\Omega}^p_X \simeq \pi_\ast\Omega^p_{\tilde{X}};$$
\item If $X=Y/G$ with $Y$ a complex manifold and $G$ a finite group of biholomorphic transformations of $Y$, 
then $$\tilde{\Omega}^p_X \simeq \left(\rho_\ast\Omega^p_Y\right)^G.$$
\end{enumerate}
\end{theorem}
See \cite[Sect.~1]{Ste77a}. Moreover we have 
\begin{theorem}
For any V-manifold $X$ the complex $\tilde{\Omega}^\bullet_X$ is a resolution of the constant sheaf $\C_X$ 
and if moreover $X$ is compact and has an algebraic structure, the Hodge spectral sequence 
\begin{equation}\label{HodgeSS}
E_1^{pq}=H^q(X,\tilde{\Omega}^p_X) \Rightarrow H^{p+q}(X,\C)
\end{equation}
degenerates at $E_1$; this gives  the Hodge filtration on the cohomology of $X$.
\end{theorem}

For spaces with arbitrary singularities the properties of these sheaves of differentials cannot 
all be preserved. Depending on the property one prefers there is a different generalisation of the above complex. 

\subsection*{The filtered de Rham complex} 
First let us focus on the Hodge-theoretic property. A filtered complex $(\tilde{\Omega}^\bullet_X,F)$ 
which resolves the constant  sheaf $\C_X$ and such that the Hodge spectral sequence (\ref{HodgeSS}) degenerates at $E_1$
has been constructed by Du Bois \cite{DB}. In the general case however the filtration 
$F$ is no longer  cutting off of the complex such that the graded complex $\tilde{\Omega}^p_X:=
\mathrm{Gr}_F^p\tilde{\Omega}^\bullet[-p]$ would be a 
single sheaf placed in degree $p$, but $\tilde{\Omega}^p_X$ is actually a complex with cohomology sheaves 
which may be non-zero on the whole range $0 \leq j \leq n-p$. This filtered complex $(\tilde{\Omega}^\bullet_X,F)$ 
is called the \emph{filtered de Rham complex} of $X$. 

Let $X'$ denote the weak normalization of $X$. It is a complex 
variety over $X$ which is homeomorphic to $X$, and sections of $\calO_{X'}$ over an open set of $X$ consist 
of those continuous functions whose restriction to the regular locus of $X$ are holomorphic. 
 
Suppose that $X$ has isolated singularities only and that $\pi:Y\to X$ is a good resolution of singularities. 
This means that the inverse image $\pi^{-1}(\Sigma)=:E$, where $\Sigma$ is the singular locus of $X$, is 
a divisor with normal crossings on $Y$. Then the cohomology sheaves of the 
filtered de Rham complex of $X$ can be described as follows: 
\begin{itemize}
\item $H^0(\tilde{\Omega}^0_X)=\calO_{X'}$;
\item $H^q(\tilde{\Omega}^p_X)=R^q\pi_\ast\Omega^p_Y(\log E)(-E)$ when $(p,q)\neq (0,0)$.
\end{itemize}
Here $\Omega^p_Y(\log E)(-E)$ is the kernel of the natural map $\Omega^p_Y \to \Omega^p_E/\mbox{torsion}$, or alternatively, 
the twist of the sheaf $\Omega^p_Y(\log E)$ of logarithmic $p$-forms with the ideal sheaf $\calO_Y(-E)$ of $E$. 
For $q>0$ the sheaf $R^q\pi_\ast\Omega^p_Y(\log E)(-E)$ has support on $\Sigma$ and its stalk at a point 
$x\in\Sigma$ has finite length $b^{p,q}$. These are the \emph{Du Bois invariants} of the isolated singularity $X,x)$. 
See \cite{Ste97}.
\subsection*{Barlet differentials} From the point of view of duality, in the hypersurface case a natural notion of
holomorphic forms turns out to be one defined as residues of a meromorphic form in the ambient space with 
logarithmic poles along the hypersurface. 
These holomorphic q-forms are the
sections of Barlet's sheaf $\omega^q_V$ , which for $q = n$ is the same as the
Grothendieck dualizing sheaf, see \cite{HP} and \cite{AT} for details.

\section{Barlet forms for hypersurfaces}
Let $P$ be a compact complex manifold of dimension $n+1$ and let $\calL$ be a line bundle on $P$. 
We consider a  hypersurface $X =V(F) \subset P$ with at most isolated singularities, given as the zero set 
of a global section $F$ of $\calL$. 
Let $\Omega^k(\ell X)$ denote the sheaf of germs of meromorphic $k$-forms on $P$ 
with poles of order at most $k$ along $X$. We have inclusions 
$\Omega^k(\ell X) \subset \Omega^k((\ell+1) X)$ and differentiation 
$d:\Omega^k(\ell X) \to \Omega^{k+1}((\ell+1) X)$. We define 
$$
\Omega^k(\log X) = \mathrm{ker}(d: \Omega^k(X) \to \Omega^{k+1}(2X)/\Omega^{k+1}(X).
$$
If $X$ is smooth then the map $d:\Omega^n(X) \to \Omega^{n+1}(2X)/\Omega^{n+1}(X)$ is surjective; 
if $X$ has isolated singularities, then the cokernel of this map is a skyscraper sheaf 
concentrated at the singular points of $X$. Its stalk at $x\in X$ is canonically isomorphic 
to $\Omega^{n+1}_X\otimes \calL^2$
which in turn is non-canonically isomorphic to the quotient of 
$\calO_{X,x}$ by the ideal generated by a local equation $f$ of $X$ and the partial derivatives of $f$. 
This stalk has finite length $\tau(X,x)$, the \emph{Tjurina number} of $(X,x)$. Then 
$$\chi(\Omega^{n+1}_X\otimes \calL^2)=\tau:=\sum_{x\in X}\tau(X,x).$$

We have the resolution 
$$
0 \to \Omega^n(\log X) \to \Omega^n(X)\stackrel{d}{\to} 
\Omega^{n+1}(2X)/\Omega^{n+1}(X) \to \Omega^{n+1}_X\otimes\calL^2 \to 0
$$
Consider the sheaf $\omega^{n-1}_X$ on $X$ defined by 
$$\omega_X^{n-1} := \Omega^n(\log X)/ \Omega^n$$
and define $$c(\calL):= \chi(\Omega^n\otimes\calL)-\chi(\Omega^n)-\chi(\Omega^{n+1}\otimes\calL^2)
+\chi(\Omega^{n+1}\otimes\calL).$$

\begin{theorem}\label{euler}
$$\chi(\omega_X^{n-1}) = \tau+c(\calL).$$
\end{theorem}
This follows immediately from the exact sequences above. 
\section{Smoothing: the specialization sequence}
We keep the notations of the previous section. Moreover we suppose that $P$ is K\"ahler. 
Let $G$ be a global section of $\calL$ which does not vanish at the singularities of $X$. 
Then there exists $\epsilon>0$ such that for $t\in \C$ with $0 <|t|<\epsilon$ the hypersurface 
$X_t$ given by $F+tG=0$ is smooth. We have the exact sequence
\begin{equation}\label{special}
0 \to H^1(X,\tilde{\Omega}^{n-1}_X) \to \mathrm{Gr}_F^{n-1}H^n(\psi\C) 
\to \mathrm{Gr}_F^{n-1}H^n(\phi \C) \to H^2(X,\tilde{\Omega}^{n-1}_X) \to 0
\end{equation}
obtained by taking $\mathrm{Gr}_F^{n-1}$ from the specialisation sequence of the smoothing (cf. \cite[Sect.~3]{Ste77a}). 
Note that it implies 
\begin{theorem}
$$\chi(\tilde{\Omega}^{n-1}_X)=s_{n-1}+c(\calL)$$
where the invariant $s_{n-1}=\sum_{x\in X}s_{n-1}(X,x)$ is a sum of local contributions from each singularity: 
 $s_{n-1}(X,x)=\dim \mathrm{Gr}_F^{n-1}H^n(\phi\C)_x$.
\end{theorem}
This theorem is analogous to Theorem \ref{euler}. 
\begin{corollary}
Suppose that $X\subset P$ is a hypersurface with isolated singularities such that 
$\mathrm{Gr}_F^{n-1}\tilde{\Omega}^\bullet_X[n-1] \simeq \omega_X^{n-1}$. Then $s_{n-1}=\tau$.
\end{corollary}
In the next section we will see that the converse of this corollary also holds. 
\subsection*{Singularity spectrum}
For an isolated hypersurface singularity $f:(\C^{n+1},0)\to (\C,0)$ we define its \emph{Milnor module}
$$\Omega_f:=\Omega^{n+1}/df\wedge\Omega^n.$$
It carries a decreasing filtration $V^\bullet$ indexed by rational numbers $a$. 
The singularity spectrum is defined in terms of the V-filtration on
$\Omega_f$ as follows: for $b\in\Q$ let $d(b):=\dim_\C\mathrm{Gr}_V^b\Omega_f$. We
put $$\mathrm{Sp}(f):=\sum_{b\in\Q}d(b)(b) \in \Z[\Q]$$
where the latter is the integral group ring of the additive group of the
rational numbers. It is called the \emph{singularity spectrum of $f$}. 

The Hodge numbers $s_k$ of the Milnor fibre of $f$ are expressed in terms of the singularity spectrum 
of $f$ by the formula 
$$s_k = \sum_{n-k-1<b\leq n-k} d(b).$$
See \cite{SS85} for details. 

\section{Local comparison}
In this section we will derive a direct relation between the Barlet and Du Bois differentials for a complete variety 
$X$  with isolated singularities. 

Let $\pi:Y \to X$ be a good resolution: if $\Sigma$ is the set of singular points of $X$, 
then $\pi^{-1}\Sigma = E$ is a divisor with normal crossings on $Y$ with smooth irreducible components. 
In this case the graded quotients of the filtered de Rham complex of $X$ are given by 
 $$\tilde{\Omega}^p_X = R\pi_\ast\Omega^p_Y(\log E)(-E)$$
for $p\geq 1$ and $\tilde{\Omega}^0_X$ is the single complex associated to the sequence 
$R\pi_\ast\calO_Y \to R\pi_\ast\calO_E \to \C_\Sigma$. By \cite{AT} we have 
$$\omega^q_X \simeq Hom(\Omega^{n-q}_X,\omega^n_X)$$
and $\omega^n_X=\omega_X$ is the dualizing sheaf. 

As we are interested in the case $p=n-1$ the first case to consider is where $X$ is a curve and $n=1$. 
Then $\tilde{\Omega}^0_X = \calO_{X'}$ where $X'$ is the \emph{weak normalization} of $X$. 
Let us compare this with the sheaf $\omega_X^0=Hom(\Omega^1_X,\omega_X)$.  

This question was first considered in the plane curve case 
by Kyoji Saito \cite{Sa80}. He proved that in this case the sheaf 
$\Omega^1(\log X)$ is locally free. 
\begin{lemma} Let $X$ be a plane curve. Then $\omega_X^{n-1}$ is torsion free of rank one. 
There is a natural injection $\tilde{\calO}_X \to \omega_X^0$, the quotient $\omega_X^0/\calO_X$ 
is concentrated in the singular locus of $X$ and its stalk at $x\in X$ has length $\tau(X,x)$.
Suppose $X=V(f)\subset (\C^2,0)$ is a reduced quasi-homogeneous plane  curve singularity: 
$w_xxf_x+w_yyf_y=f$. Then the forms $\frac{df}{f}$ and $\frac{w_yydx-w_xxdy}{f}$ form a local basis of 
$\Omega^1(\log X)$. If $f$ is not quasi-homogeneous, then $\frac{df}{f}\in (x,y)\Omega^1(\log X)$.
\end{lemma}
{\em Proof } On the regular locus of $X$ we have an isomorphism between $\omega_X^0$ and $\calO_X$ by the residue map. 
Let us check that it extends to the desired injection. 
This, and the fact that $\omega_X^0$ is torsion free, can be checked locally near every singular point. 
So consider the case of a reduced plane curve singularity $(X,0) \subset (\C^2,0)$ 
given by a squarefree function germ $f\in\C\{x,y\}$. 
We have $\C\{x,y\}dx \oplus \C\{x,y\}dy \subset \Omega^1(\log X)_0 \subset \C\{x,y\}dx/f \oplus \C\{x,y\}dy/f$ 
so writing $\calO = \calO_{X,0}$ we have $\omega^0_{X,0} \subset \calO dx/f \oplus \calO dy/f$, 
so it is a subsheaf of a locally free sheaf. 
Hence $\omega_X^0$ is torsion free, and its rank is the same as the rank of its restriction to the regular locus, 
which equals one. 

The element $a\ dx/f + b\ dy/f$ with $a,b\in\calO$ belongs to $\omega_{X,0}^0$ iff $a\eta-b\xi=0$ 
where $\xi,\eta$ are the images of $f_x,f_y$ in $\calO$ respectively. 
This equation has the obvious solution $a=\xi,b=\eta$ 
which corresponds to the germ $df/f$ and hence to $1\in\calO$. 
The injection $\calO \hookrightarrow \omega_X^0$ is therefore given by $c \mapsto c\xi\ dx/f + c\eta\ dy/f$. 

To compute the length of the stalk of $\omega_X^0/\calO_X$ at a singular point of $X$ we use a global argument, 
even if the question is local. To this end, suppose that $x\in X$ 
is the unique singular point of a plane projective curve of degree $d$. 
Then the length of   $\omega_X^0/\calO_X$ at $x$ is the same as the Euler Poincar\'e characteristic 
$\chi(\omega_X^0/\calO_X) = \chi(\omega_X^0) - \chi(\calO_X)$.  
Recall that $\chi(\omega_X^0) = \tau+1 -c_d$ whereas $\chi(\calO_X)=1-c_d$. Hence $\chi(\omega_X^0/\calO_X)=\tau$. 

To show that $\omega_X^0$ contains $\tilde{\calO}_X$ 
observe that for any $\omega\in\Omega^1_X$ and $a\in\tilde{\calO}_X$, 
the product $a\omega$ lies in $\Omega_{\tilde{X}}$ so has no residues; 
hence it belongs to $\omega_X$, i.e. $a \in Hom(\Omega^1_X,\omega_X)$. 

The remaining statements are left to the reader. See also \cite[Proof of Theorem 2.11]{Sa80}. 
QED

\vspace{4mm}\noindent
{\bf Example } Let $X$ be a projective plane curve with only ordinary double points. 
Then $\omega_X^0 \simeq \tilde{\calO}_X$.
In particular $H^0(X,\calO_X) \to H^0(X,\omega_X^0)$ is an isomorphism if and only if $X$ is irreducible. 

Recall that $X'$ denotes the weak normalization of the curve $X$. 

\begin{theorem} \label{noplanecurvesing}
Suppose that $X$ is a plane curve such that $\calO_{X'} = \omega_X^0$. Then $X$ is smooth. 
\end{theorem}
{\em Proof } If  $\calO_{X'} = \omega_X^0$ then $\calO_{X'}/\calO_X$ and $\omega_X^0/\calO_X$ 
have stalks of the same lengths at all singular points. 
Hence for such a singularity one has the equality $\delta-r+1=\tau$.  
By \cite[Lemma 6.1.2 and Cor.~6.1.4]{BuGr} $\tau\geq \delta+m-r$ where $m$ is the multiplicity. 
Hence $m=1$ so $X$ has no singular point. 

Here is another argument, based on the spectrum. 
Consider $Q^f = \Omega^2/df\wedge\Omega^1$ with its spectral $V$-filtration. 
We have $fQ^f \subset V^{>0}Q^f$ so 
$\tau = \dim Q^f/fQ^f \geq \dim Q^f/V^{>0} = \delta$ 
with equality iff $r=1$ and $fQ^f=V^{>1}$. 
So $X$ is an irreducible plane curve singularity, 
and multiplication by $f$ gives an isomorphism $Q^f/V^{>0}  \to V^{>0}$. 
If $\alpha_1,\ldots,\alpha_{\delta}$ are the positive spectral numbers of $f$ in increasing order, 
then the spectrum of $f$ is $-\alpha_\delta,\ldots,-\alpha_1,\alpha_1,\ldots,\alpha_\delta$. 
We find that $\alpha_j +\alpha_{\delta-j+1} \geq 1$ for all $j$. 
This implies that the surface singularity with equation $f(x,y)+z^2=0$ 
has geometric genus $p_g \geq \mu/4$, 
but N\'emethi \cite{N} has shown that for such a surface singularity $p_g\leq \mu/6$. 
This means that the spectral numbers have to lie closer to the middle than forced by the condition $fQ^f = V^{>0}Q^f$. 

Another argument is based on Hertling's conjecture \cite{H} on the variance of the spectrum, 
which has been proved by Br\'elivet in the curve case \cite{B}. 
If $\alpha_j+\alpha_{\delta-j+1} \geq 1$, then $\alpha_j^2+\alpha_{\delta-j+1}^2 \geq \frac 12$ 
so $\sum_{i=1}^\mu \alpha_j^2 \geq \mu/4$. On the other hand, by Hertling's conjecture 
$$\frac{1}{\mu}\sum_{i=1}^\mu \alpha_j^2 \leq \frac{1}{12}(\alpha_{\mu}-\alpha_1)\leq \frac 16.$$
QED

\vspace{4mm}\noindent
Next we turn to the study of $\omega_X^1$ in the case $n\geq 2$. 
Note that $\omega_X^1$ coincides with $\Omega_X^1$ on $X\setminus \Sigma$. 
Moreover $\omega_X^1$ fits in the exact sequence 
$$0 \to \omega_X^1 \to \Omega^n(X)\otimes \calO_X \to \Omega^{n+1}(2X)\otimes \calO_X $$
hence $\omega_X^1 \simeq j_\ast\Omega_X^{n-1}$ where $j:X\setminus\Sigma \hookrightarrow X$. 
We see that $\omega_X^1 \simeq \tilde{\Omega}_X^{n-1}$ if $X$ is a V-manifold. 
Let us look for the class of singularities which one may admit for this to be true. 
Consider a good resolution $\pi:(Y,E) \to (X,x)$ of an isolated $n$-dimensional singularity. 
Define 
$$q'(X,x) = \ell\left(j_\ast\Omega_X^{n-1}/\pi_\ast\Omega^{n-1}_Y(\log E)(-E)\right)_x$$
and the Du Bois invariant (cf. \cite{Ste97})
$$b^{n-1,1}(X,x) = \ell\left(R^1\pi_\ast\Omega^{n-1}_Y(\log E)(-E)\right)_x.$$
Then clearly one has the 
\begin{theorem} Let $X$ be an $n$-dimensional complex space with only isolated singularities, with $n\geq 2$. 
The following are equivalent:
\begin{enumerate}
\item $\omega_X^{n-1} = \tilde{\Omega}_X^{n-1}$;
\item $q'(X,x) = b^{n-1,1}(X,x) = 0$ for each singular point of $X$.
\end{enumerate}
\end{theorem}
Indeed, if $q'=0$ then $$\omega_X^{n-1}=j_\ast\Omega_X^{n-1} = \pi_\ast\Omega^{n-1}_Y(\log E)(-E)$$ 
and if moreover $b^{n-1,1}=0$ then $R^1\pi_\ast\Omega^{n-1}_Y(\log E)(-E)=0$ so 
$\pi_\ast\Omega^{n-1}_Y(\log E)(-E) = R\pi_\ast\Omega^{n-1}_Y(\log E)(-E)$. Conversely, the equality 
$\omega_X^{n-1} = \tilde{\Omega}_X^{n-1}$ implies equality of their Euler characteristics, whose difference is 
equal to $q'+b^{n-1,1}$. 
\begin{corollary}
Suppose that $(X,x)$ is an isolated hypersurface singularity. Then 
$$q'(X,x)+b^{n-1,1}(X,x)= \tau(X,x)-s_{n-1}(X,x).$$
\end{corollary}

Next we investigate which surface singularities have $q'(X,x) = b^{1,1}(X,x) = 0$. 
First recall the following result of Wahl \cite[Corollary 2.9]{W85} :
\begin{theorem}
For a two-dimensional smoothable normal Gorenstein singularity with Milnor fibre $F$ 
write $\mu=\mu_0+\mu_+ +\mu_- $ from diagonalizing the intersection pairing on $H_2(F,\R)$. Then 
$$\tau\geq \mu_0+\mu_-=\mu-(2p_g-2g-b).$$
\end{theorem}\label{thm:Wahl}
Here $p_g$ is the geometric genus, $b$ is the first Betti number of the dual graph of a good resolution 
and $g$ is the sum of the genera of the irreducible components of its exceptional divisor. 
\begin{theorem}
For a two-dimensional smoothable normal Gorenstein singularity the following are equivalent:
\begin{enumerate}
\item $b^{1,1}=q'=0$;
\item $\tau = \mu_0+\mu_-$ and $g=0$.
\end{enumerate}
\end{theorem}
{\em Proof } 
In the surface case we have $\mu=s_0+s_1+s_2$ where $s_2=p_g$ and $s_0=p_g-g-b$. If 
$b^{1,1}=q'=0$ then $\tau=s_1=\mu-(2p_g-g-b)=\mu_0+\mu_--g$ hence by Theorem \ref{thm:Wahl} we have 
$g=0$ and $\tau=\mu_0+\mu_-$. 
The converse implication is similar.

\begin{corollary}
The following surface singularities satisfy $q'=b^{1,1}=0$:
\begin{enumerate}
\item rational double points (ADE-singularities)
\item cusps (singularities of type $T_{pqr}$ with $\frac 1p +\frac 1q +\frac 1r <1$);
\item generic $\mu$-constant deformations of $z^2+x^{2a+1}+y^{2a+2}$ (those which have minimal Tjurina number).
\end{enumerate}
\end{corollary}
Indeed, the rational double points have $\mu=\tau=\mu_-$, whereas the cusps have $\mu_+=1=\mu-\tau$. 
Finally the last category of examples was considered in \cite[Example 4.6]{W85} 
and shown by Zariski to have $\tau=3a(a+1)=\mu-2p_g$. 
QED

\vspace{4mm}\noindent {\bf Example }
According to SINGULAR, the singularity $x^7+x^4y^2+x^2y^4+y^7+z^2$ has $\mu=27$ and $\tau=23$. 
Moreover, $p_g=3$ and $b=2,\ g=0$, and 
$\mu-(2p_g-b)=27-4=\tau$ so $q'=b^{1,1}=0$.

\vspace{4mm}\noindent {\bf Example }
The singularity $x^3+y^{10}+z^{19}$ is considered in \cite[Sect.~5]{LP}. 
A generic $\mu$-constant deformation has $\tau=246$ whereas $\mu=324$ and $p_g=39$, $g=b=0$. 
so $\mu-\tau=2p_g$ and $q'=b^{1,1}=0$.

\vspace{4mm}\noindent {\bf Example }
It is not always so that for a generic $\mu$-constant deformation of a quasi-homogeneous surface sngularity 
one has $q'=b^{1,1}=0$. The exceptional unimodal non-quasi-homogeneous singularities have 
$\tau=\mu-1,\ p_g=1$ and $g=b=0$. So $\mu_0=0,\ \mu_-=\mu-2$. 

Take  the singularity $x^5+y^{11}+z^2$, also considered in \cite[Sect.~5]{LP}. 
It has $\mu=40$ and $\tau_\mathrm{min}=34$ whereas $p_g=4,\ g=0$. 
That $\mu-\tau\leq 6$ can be seen from the spectral numbers. 
The submodule $fQ^f$ of $Q^f=\Omega^3/df\wedge\Omega^2$ is cyclic with generator 
$[f\omega]\in V^{>\frac{87}{110}}=V^{\frac{89}{110}}$. But then 
$$[xf\omega] \in V^{\frac{111}{110}},\, [yf\omega] \in 
V^{\frac{99}{110}}= V^{\frac{101}{110}},\, [y^2f\omega] \in V^{\frac{111}{110}}$$
so the spectral numbers of the filtration of $fQ^f$ induced by $V$ have the gaps $\frac{91}{111}$ and $\frac{93}{110}$. 

\vspace{4mm}\noindent {\bf Question } Is this kind of lower bound for $\tau$ on the 
$\mu$-constant stratum provided by the spectral numbers sharp?

\vspace{4mm}\noindent {\bf Remark }
The only example of an isolated hypersurface singularity in dimension $n\geq 3$ 
I know that satisfies $q'=b^{n-1,1}=0$ is the ordinary double point in dimension three. 
If Hertling's conjecture is valid, I can prove that no such singularities exist in dimension 
$\geq 9$, and no rational singularity in dimension $\geq 6$. 

Consider the special case of ``double suspension'' singularities $g=f(x,y)+zw$ with $f$ squarefree. 
These are rational, and they belong to our class iff $\tau_f=\delta_f$. 
By the inequality $\tau\geq \delta+m-r$ for curve singularities this implies that $r=m$ 
so $f$ is a $\mu$-constant deformation of a homogeneous singularity of degree $m$. 
This has finite order monodromy and highest spectral number $1-\frac2m$ 
and this implies that $V^{-\frac 2m}Q^f$ is in the kernel of multiplication by $f$. 
Hence $\tau_f\geq \dim V^{-\frac 2m}Q^f = \delta_f +2m-5 $ so $\tau_f-\delta_f\geq 1$ unless $m=2$ 
and we have the ordinary double point! 
\section{Application to projective hypersurfaces}
In this section we come back to the problem mentioned in the introduction: investigate the relation between 
the cohomology of a projective hypersurface with isolated singularities and certain graded parts of its Jacobian ring. 

We consider a hypersurface $X =V(F) \subset \PP^{n+1}$ of degree $d$ with at most isolated singularities. 
Recall that $\Omega^k(\ell X)$ is the sheaf of germs of meromorphic $k$-forms on $\PP^{n+1}$ with poles of 
order at most $k$ along $X$. We have inclusions 
$\Omega^k(\ell X) \subset \Omega^k((\ell+1) X)$ and differentiation $d:\Omega^k(\ell X) \to 
\Omega^{k+1}((\ell+1) X)$. We defined 
$$
\Omega^k(\log X) = \mathrm{ker}(d: \Omega^k(X) \to \Omega^{k+1}(2X)/\Omega^{k+1}(X).
$$
If $X$ is smooth then the map $d:\Omega^n(X) \to \Omega^{n+1}(2X)/\Omega^{n+1}(X)$ is surjective; 
if $X$ has isolated singularities, then the cokernel $\Omega^{n+1}_X(2X)$ of this map is a skyscraper sheaf concentrated 
at the singular points of $X$. Its stalk at $x\in X$ is isomorphic to the quotient of $\calO_{X,x}$ 
by the ideal generated by a local equation $f$ of $X$ and the partial derivatives of $f$. This stalk 
has finite length $\tau(X,x)$, the \emph{Tjurina number} of $(X,x)$. 

By Bott's vanishing theorem \cite{Bo57} $H^i(\Omega^n(X))=H^i(\Omega^{n+1}(2X)/\Omega^{n+1}(X))=0$ for 
$i>0$, so we have a resolution 
$$
0 \to \Omega^n(\log X) \to \Omega^n(X)\stackrel{d}{\to} \Omega^{n+1}(2X)/\Omega^{n+1}(X) \to \Omega^{n+1}_X(2X) \to 0
$$
of $\Omega^n(\log X)$ by sheaves which are acyclic, hence the cohomology groups of the complex of 
global sections of these sheaves 
$$ 0 \to H^0(\Omega^n(X)) \to H^0(\Omega^{n+1}(2X)/\Omega^{n+1}(X)) \to H^0(\Omega^{n+1}_X(2X)) \to 0$$
are isomorphic to the cohomology groups of $ \Omega^n(\log X)$. Explicitly we have the complex 
$$
0 \to \calS_{d-n-2} \stackrel{E}{\to} \calS_{d-n-1}^{\oplus n+2} \stackrel{h}{\to} \calS_{2d-n-2}/F\calS_{d-n-2} 
\to H^0(\Omega^{n+1}_X(2X)) \to 0
$$
where $E(B)=(X_0B,\ldots,X_{n+1}B)$ (corresponding to the Euler vector field) and 
$$h(A_0,\ldots,A_{n+1})= \sum_{i=0}^{n+1}A_i\frac{\partial F}{\partial X_i} \ \mathrm{mod}\ F
$$ 

If $X$ is smooth, we have the residue exact sequence 
$$
0 \to \Omega^n \to \Omega^n(\log X) \to \Omega^{n-1}_X \to 0
$$
by which the cohomology groups of $ \Omega^n(\log X)$ are identified with the primitive 
cohomology groups  $H^{n-1,1}_\mathrm{prim}(X)$. 
\begin{lemma}
If $X$ is smooth, then 
$$\sum_{n\geq 0}\dim \calR_nt^n=\left(\frac{t^{d-1}-1}{t-1}\right)^{n+2}.$$
\end{lemma}
The proof uses the fact that the partials of $F$ form a regular sequence in $\calS$, so 
we have the Koszul complex resolving $\calR$. 

\begin{proposition}
Let $\calR = \calS/\calJ(F)$ as in the introduction. Then 
$$
\dim \calR_{2d-n-2} = c_d - \dim H^0(X,\Omega^n(\log X)).
$$ 
Moreover, the map $H^i(\Omega^n(log X)) \to H^i(\omega^{n-1}_X)$ is an isomorphism for $i\neq n-1$ and we have the exact sequence 
$$ 0 \to H^{n-1}(\Omega^n(\log X)) \to H^{n-1}(\omega^{n-1}_X) \to H^n(\Omega^n) \to 0.$$
\end{proposition}
{\em Proof } Note that $$\calR_{2d-n-2} = \mathrm{coker}(\calS_{d-n-1}^{\oplus n+2} \stackrel{h}{\to} 
\calS_{2d-n-2}/F\calS_{d-n-2})$$
and that $\mathrm{ker}(h)=\mathrm{im}(E)$ in the smooth case. Hence 
$$c_d = \dim  \calS_{2d-n-2}- (n+2)\dim \calS_{d-n-1}= {2d-1 \choose n+1} - (n+2){d \choose n+1}$$
and 
$$\mathrm{ker}(h)/\mathrm{im}(E) =  H^0(X,\Omega^n(\log X))) = \mathrm{ker}\left(\calS_{d-n-1}^{\oplus(n+2)} 
\to \calS_{2d-n-2}\right).$$ To prove the remaining statements, note that $H^i(\Omega^n)=0$ for all $i\neq n$. 
So we only have to show that 
$$H^n(\Omega^n) \to H^n(\Omega^n(\log X))$$
is the zero map. We do this by induction on $n$. The case $n=1$ is obvious. Let $n\geq 2$. 
Consider a general hypersurface $L \subset \PP^{n+1}$; we have the commutative diagram with exact rows 
$$
\begin{array}{ccccccccc}
0 & \to & \Omega^n & \to & \Omega^n(\log L) & \to & \Omega_L^{n-1} & \to & 0 \\
&& \downarrow && \downarrow && \downarrow &&\\
0 & \to & \Omega^n(\log X) & \to & \Omega^n(\log L+X) & \to & \Omega_L^{n-1}(\log X\cap L) & \to & 0 
\end{array}
$$
which gives rise to the commutative diagram 
$$
\begin{array}{ccc}
H^{n-1}(\Omega_L^{n-1}) & \stackrel{a}{\rightarrow} & H^{n-1}(\Omega_L^{n-1}(\log X\cap L)) \\
\downarrow b && \downarrow \\
H^n(\Omega^n) & \stackrel{c}{\rightarrow} & H^n(\Omega^n(\log X)) 
\end{array}
$$
By Lefschetz' theory, the Gysin map $b$ is an isomorphism, and by induction hypothesis $a$ is the zero map. 
Hence $c$ is the zero map. 
QED
\begin{corollary}
Suppose that $H^0(\Omega^n(\log X))=0$. Then $\calR_{2d-n-2}$ has the expected dimension $c_d$ and we have the exact sequence 
\begin{equation}\label{Fseq}
0 \to H^1(\Omega^n(\log X)) \to \calR_{2d-n-2} \to H^0(X,\Omega^{n+1}_X(2X)) \to H^2(\Omega^n(\log X))\to 0
\end{equation}
Moreover $H^i(\Omega^n(\log X))=0$ for all $i\geq 2$. 
\end{corollary}
From now on we suppose that $\omega^{n-1}_X \simeq \tilde{\Omega}^{n-1}_X$, i.e. $q'=b^{n-1,1}=0$ for all 
singular points of $X$. Moreover we will suppose that $n\geq 2$, as in the case $n=1$ there are no singular points on $X$. 
This guarantees that we have isomorphisms 
\begin{equation}\label{equalities}
H^i( \Omega^n(\log X)) \simeq \mathrm{Gr}_F^{n-1}H^{n-1+i}(X)_\mathrm{prim}.
\end{equation}

\begin{lemma}
$H^i(\PP^{n+1},\Omega^n(\log X))=0$ for all $i\neq 1,2$ 
\end{lemma}
{\em Proof } This follows from (\ref{equalities}) and the fact that for a hypersurface $X$
in $\PP^{n+1}$ with isolated singularities one has $H^k(\PP^{n+1},\Q)\simeq H^k(X,\Q)$ for all $k\neq n,n+1,2n$. 
\begin{corollary} $\dim\calR_{2d-n-2}=c_d$. Moreover we have the exact sequence 
\begin{equation}\label{Fseq}
0 \to H^1(\Omega^n(\log X)) \to \calR_{2d-n-2} \to H^0(X,\Omega^{n+1}_X(2X)) \to H^2(\Omega^n(\log X))\to 0
\end{equation}
\end{corollary}
Let $G$ be a homogeneous form of degree $d$ which does not vanish at the singularities of $X$. 
Then there exists $\epsilon>0$ such that for $t\in \C$ with $0 <|t|<\epsilon$ the hypersurface 
$X_t$ given by $F+tG=0$ is smooth. We let $\calR^t$ denote its Jacobian ring. 
Note that $\lim_{t\to 0}\calJ(F+tG)_k$ makes sense in the Grassmannian of $\calS_k$, 
and that it contains $\calJ(F)_k$. Hence, as $\calR_k$ and $\calR^t_k$ have equal dimension, they are equal. 
So we have the exact sequence
\begin{equation}\label{special}
0 \to H^1(X,\tilde{\Omega}^{n-1}_X)_{\mathrm{prim}} \to \calR_{2d-n-2} \to 
\mathrm{Gr}_F^{n-1}H^n(\phi C) \to H^2(X,\tilde{\Omega}^{n-1}_X)_{\mathrm{prim}} \to 0
\end{equation}
Under our hypotheses, the sequences (\ref{Fseq}) and (\ref{special}) are identical! 

\vspace{4mm}\noindent
{\bf Remark } Our reasoning also applies to hypersurfaces in weighted projective spaces, as long as they are transverse 
to the singular strata (so have isolated singularities only at regular points of the ambient space).


\begin{thebibliography}{999999}
\bibitem{AT} A.G.~Aleksandrov, A.T.~Tsikh, Th\'eorie des r\'esidus de Leray et formes de Barlet 
sur une intersection compl\`ete singuli\`ere, {\sl C.R.~Acad.~Sci.~Paris} {\bf 333} S\'erie I (2001), 973--978.

\bibitem{Bo57} R.~Bott, Homogeneous vector bundles, {\sl Annals of Math.} {\bf 66} (1957), 203--248. 

\bibitem{B} Th.~Br\'elivet, The Hertling conjecture in dimension two. math.AG/0405489.  

\bibitem{Br70} E.~Brieskorn, Singular elements of semi-simple algebraic groups. 
Proc.~Internat.~Congr.~Math., Nice 1970, Vol.~2, Gauthier-Villars, Paris, 1971, pp.~279--284.

\bibitem{BuGr} R.-O.~Buchweitz, G.-M.~Greuel, The Milnor number and deformations of 
complex curve singularities. {\sl Inventiones math.} {\bf 58} (1980), 241--281.

\bibitem{DB} Ph.~du Bois, Complexe de de Rham filtr\'e d'une vari\'et\'e singuli\`ere. 
{\sl Bull.~Soc.~Math.~France} {\bf 109} (1981), 41--81.


\bibitem{Gre80} G.-M.~Greuel, Dualit\"at in der lokalen Kohomologie isolierter Singularit\"aten, 
{\sl Math.~Ann.} {\bf 250} (1980), 157-173.

\bibitem{Gr69} Ph.~A.~Griffiths, On the periods of certain rational integrals: I and II. 
{\sl Annals of Math.} {\bf 90} (1969), 460--495 and 498--541. 

\bibitem{HP} G.~Henkin, M.~Passare, 
Abelian differentials on singular varieties
and variations on a theorem of Lie-Griffiths, {\sl Invent. math.} {\bf 135} (1999), 297--328.

\bibitem{H} C.~Hertling, Frobenius manifolds and the variance of the spectral numbers, 
In: {\sl New Developments in Singularity Theory}, D.~Siersma, C.T.C.~Wall and V.~Zakalyukin eds. 
Kluwer Academic Publishers 2001.

\bibitem{LP} O.A.~Laudal, G.~Pfister, Local Moduli and Singularities, Lecture Notes in Math.~ 1310, 
Springer-Verlag, Berlin etc. 1988. 

\bibitem{N} A.~N\'emethi, Dedekind sums and the signature of $f(x,y)+z^N$, II, 
{\sl Selecta Mathematica, New Series} {\bf 5} (1999), 161-179. 

\bibitem{Sa80} K.~Saito, Theory of logarithmic differential forms and logarithmic vector fields, 
{\sl J.~Fac.~Sci.~Univ.~ Tokyo, Ser.~IA} {\bf 27}(2) (1980), 265--291.

\bibitem{SS85} J.~Scherk, J.H.M.~Steenbrink, On the Mixed Hodge Structure on the Cohomology of the
Milnor Fibre, {\sl Math.~Annalen} {\bf 271} (1985), 641--665.

\bibitem{Ste77a} J.H.M.~Steenbrink, Mixed Hodge structures on the
vanishing cohomology. In {\sl Real and Complex Singularities,
Oslo, 1976}, Sijthoff-Noordhoff, Alphen a/d Rijn, 525--563 (1977).

\bibitem{Ste97} J.H.M.~Steenbrink, Du Bois invariants of isolated complete intersection singularities. 
{\sl Ann.~Inst.~Fourier, Grenoble} {\bf 47}, 5 (1997), 1367--1377.

\bibitem{StrS85} D.~van Straten and J.H.M.~Steenbrink, Extendability of holomorphic differential forms 
near isolated hypersurface singularities. {\sl Abh.~Math.~Sem.~Univ.~Hamburg} {\bf 55} (1985), 97--110.

\bibitem{W85} J.~Wahl, A characterization of quasi-homogeneous Gorenstein surface singularities. 
{\sl Compos.~Math.} {\bf 55} (1985), 269--288. 

\end{thebibliography}
\end{document}